\documentclass[J,12pt]{amsart}
\usepackage{epsfig,epsf,amsmath,amssymb}
\newtheorem{theorem}{Theorem}[section]

\newtheorem{definition}{Definition}
\newtheorem{example}[theorem]{Example}
\newtheorem{proposition}[theorem]{Proposition}
\newtheorem{remark}[theorem]{Remark}

\begin{document}
\title[]{Digital $k$-continuity, digital $k$-isomorphism, local $k$-isomorphism, radius $2$-local $k$-isomorphism, and digital $k$-homotopy}

\author[]{Sang-Eon Han}
\address[]{Department of Mathematics Education, Institute of Pure and Applied Mathematics\\
	Jeonbuk National University, Jeonju-City Jeonbuk, 54896, Republic of Korea\\
	e-mail address:sehan@jbnu.ac.kr, Tel: 82-63-270-4449.}

\begin{abstract}
					The present paper refers to the notions of digital continuity, digital $k$-isomorphism, local $k$-isomorphism, radius $2$-local $k$-isomorphism, and digital $k$-homotopy motivated by the Khalimsky's version.
	We discuss something incorrectly mentioned in Boxer's papers and suggest some 
		accurate information. 
\end{abstract}


\maketitle
\newpage

\section{\bf Essential concepts for studying digital images}\label{s1}

The present paper will often follow the notation.
${\mathbb N}$ indicates the set of positive integers, ${\mathbb Z}$ means the set of  integers, $\lq\lq$$:=$" is used for introducing a new term and $\lq\lq$$\subset$" stands for a subset as usual. 
Besides,  for $a, b \in {\mathbb Z}$, 
$[a, b]_{\mathbb Z}:=\{t \in {\mathbb Z}\,\vert \,a \leq t \leq b\}$.
We also abbreviate the term $\lq\lq$digital topological" as $DT$.\\

To study digital images from the viewpoint of digital topology, we have
considered digital $k$-connectivities for the $n$-dimensional lattice set ${\mathbb Z}^n, n \in {\mathbb N}$.
In detail, before 2004 (see \cite{H1,H2,H3}), focusing on the $3$-dimensional real world, all papers in digital topology dealt with only  low dimensional digital image $(X, k)$, i.e., $X \subset {\mathbb Z}^3$ with one of the $k$-adjacency of
${\mathbb Z}^3$ \cite{B2,KR1,R1,R2}. \\

After generalizing these $k$-connectivities into those for $n$-dimensional cases, $n \in {\mathbb N}$, since 2004 (\cite{H1,H2,H3}), 
a digital image $(X, k)$ has often assumed to be a set $X$ in ${\mathbb Z}^n$ with one of the $k$-adjacency relations of
${\mathbb Z}^n$  from (2.2) below. 
Indeed, the $k$-adjacency relations for $X \subset {\mathbb Z}^n, n \in {\mathbb N}$, were initially established in \cite{H3} (see also \cite{H1,H2,H6,H11}).
 More precisely, for a natural number $t$, $1 \leq t \leq n$, the distinct points
$p = (p_i)_{i \in [1,n]_{\mathbb Z}}$  and $q=(q_i)_{i \in [1,n]_{\mathbb Z}}\in {\mathbb Z}^n$
are $k(t, n)$-adjacent (see the page 74 of \cite{H3}) if 
$$\text{at most}\,\,t\,\,\text{of their coordinates  differ by}\,\,\pm1\,\,\text{and the others~coincide.}\eqno(2.1)$$

According to this criterion, the digital $k(t, n)$-connectivities of ${\mathbb Z}^n, n \in {\mathbb N}$, are formulated \cite{H3} (see also \cite{H11}) as follows:
$$k:=k(t,n)=\sum_{i=1}^{t} 2^{i}C_{i}^{n}, \text{where}\,\, C_i ^n:= {n!\over (n-i)!\ i!}. \eqno(2.2)$$
For instance, in ${\mathbb Z}^4$, 
$k(1,4)=8$, $k(2,4)=32$, $k(3,4)=64$, and $k(4,4)=80$; in ${\mathbb Z}^5$, 
$k(1,5)=10$, $k(2,5)=50$, $k(3,5)=130$, $k(4,5)=210$, and $k(5,5)=242$; and 
in ${\mathbb Z}^6$, 
$k(1,6)=12$, $k(2,6)=72$, $k(3,6)=232$, $k(4,6)=472$, $k(5,6)=664$, $k(6,6)=728$ \cite{H1,H3,H11}.\\
The notation $\lq\lq$$k(t, n)$" of (2.2) indeed gives us essential information 
to study digital images such as the dimension $n$ and the number $t$ associated with the statement of (2.1).
 Han's approach using $k:=k(t,n)$ can be helpful to study  more complicated cases such as digital products with a $C$-compatible, a normal, and a pseudo-normal $k$-adjacency \cite{H3,H10,H18,KHL1} and their applications.\\
 
Hereinafter, a relation set $(X, k)$ is assumed in ${\mathbb Z}^n$ with a $k$-adjacency of (2.2), called a digital image. Namely, we now abbreviate $\lq\lq$digital image $(X, k)$" as $\lq\lq$$(X, k)$" if there is no confusion.\\

Let us now recall a $\lq\lq$digital $k$-neighborhood of the point $x_0$" in $(X, k)$.
For $(X, k)$ and a point $x_0 \in X$, we have often used the notation
$$ N_k(x_0, 1):=\{x \in X \,\vert\,\,x\,\,\text{is}\,\,k\text{-adjacent to}\,\,x_0\} \cup \{x_0\}, \eqno(2.3)$$
which is called a digital $k$-neighborhood of the point $x_0$ in $(X, k)$.
The notion of (2.3) indeed facilitates some studies of digital images, e.g., digital covering spaces \cite{H3,H5,H6},
 digital $(k_0,k_1)$-homotopy \cite{B3,BK1,H3,H5,}, digital-topological ($DT$-, for brevity) $k$-group \cite{H15,H16}, and so on. 

 The notation of (2.3) indeed comes from the $\lq\lq$digital $k$-neighborhood of $x_0 \in X$ with
 radius $\varepsilon \in {\mathbb N}$" in \cite{H1,H2,H3} (see (2.4) below). 
 More precisely, for $(X,k), X \subset {\mathbb Z}^n$ and two distinct points $x$ and $x^\prime$ in $X$, a {\it simple $k$-path} from $x$ to $x^\prime$ with $l$ elements is said to be a sequence $(x_i)_{i \in [0, l]_{\mathbb Z}}$ in $X$ such that $x_0=x$ and $x_l=x^\prime$ and further, $x_i$ and $x_j$ are
$k$-adjacent if and only if $\vert i-j\vert=1$ \cite{KR1}. We say that a {\it length} of the simple $k$-path, denoted by
$l_k(x, y)$, is the number $l$ \cite{H3}.

\begin{definition} \cite{H1} (see also \cite{H3}) For a digital image $(X, k)$ on ${\mathbb Z}^n$,
	the digital $k$-neighborhood of $x_0 \in X$ with
	radius $\varepsilon$ is defined in $X$ to be the following subset of $X$
	$$N_k(x_0, \varepsilon) = \{x \in X \,\vert\, l_k(x_0, x) \leq
	\varepsilon\}\cup\{x_0\}, \eqno(2.4) $$
	where $l_k(x_0, x)$ is the length of a shortest simple $k$-path from $x_0$ to $x$ and  $\varepsilon\in {\mathbb N}$.
\end{definition}

Owing to Definition 1, it is clear that $N_k(x_0, 1)$ of (2.3) is a special case of $N_k(x_0, \varepsilon)$ of (2.4).
Note that the neighborhood $N_k(x_0, 2)$ plays an important role in establishing 
the notions of a radius $2$-local $(k_1,k_2)$-isomorphism and a radius $2$-$(k_1,k_2)$-covering map which is strongly associated with the
homotopy lifting theorem \cite{H2}  and is essential for calculating digital fundamental groups of digital images $(X,k)$ \cite{BK1,H3}. \\

We say that $(Y, k)$ is $k$-connected \cite{KR1} if for arbitrary distinct points $x, y \in Y$ there is a finite sequence $(y_i)_{i \in [0, l]_{\mathbb Z}}$ in $Y$ such that
$y_0=x$, $y_l=y$ and $y_i$ and $y_j$ are
$k$-adjacent if $\vert\, i-j \,\vert=1$ \cite{KR1}.
Besides, a singleton is assumed to be $k$-connected.\\

A simple closed $k$-curve (or $k$-cycle)
with $l$ elements in ${\mathbb Z}^n, n \geq 2$, denoted by $SC_k^{n,l}$ \cite{H3,KR1}, $l(\geq 4) \in {\mathbb N}$, is defined as the sequence $(y_i)_{i \in [0, l-1]_{\mathbb Z}}\subset {\mathbb Z}^n$ such that
$y_i$ and $y_j$ are $k$-adjacent if and only if  $\vert\, i-j \,\vert=\pm1(mod\,l)$.
Indeed, the number $l(\geq 4)$ is very important in studying $k$-connected digital images 
from the viewpoint of digital homotopy theory.
The number $l$ of $SC_k^{n,l}$ can be an even or an odd number. Namely, the number $l$ depends on the dimension $n$ and the $k$-adjacencies (see Remark 1.1 below).

\begin{remark}  (see (5) of \cite{H13}) The number $l$ of $SC_k^{n,l}$ depends on the dimension $n\in {\mathbb N}\setminus \{1\}$ and the number $l \in \{2a, 2a+1\}, a \in {\mathbb N} \setminus \{1\}$.
\end{remark}

\section{\bf Digital $k$-continuity, digital $k$-isomorphism, local $k$-isomorphism, and radius $2$-local $k$-isomorphism}\label{s4}

Rosenfeld \cite{R1} initially introduced the notion of digital $k$-continuity.
More precisely, Rosenfeld defined that 
a map $f: (X, k_1) \to (X, k_2)$ is said to be $(k_1, k_2)$-continuous if  the image $(f(A), k_2)$ is $k_2$-connected for each  $k_1$-connected subset $A$ of $(X, k_1)$, where $X \subset {\mathbb Z}^3$.
Also, Rosenfeld represented it as follows:
 A function $f: (X,k_1) \to (Y, k_2)$ is $(k_1, k_2)$-continuous if and only if  two points $x$ and $x^\prime$ which are $k_1$-adjacent in $(X,k_1)$ implies that either $f(x)=f(x^\prime)$ or $f(x)$ is $k_2$-adjacent to $f(x^\prime)$, where $X$ and $Y$ are subsets of ${\mathbb Z}^3$.
 Based on this approach, Han \cite{H3} represented the $(k_1, k_2)$-continuity of a map $f: (X, k_1) \to (Y, k_2)$ with the following statement which has strongly used to study digital covering spaces and pseudo-convering spaces, $DT$-$k$-group theory \cite{H15,H16}, and so forth, where $X \subset {\mathbb Z}^{n_1}$ and $Y \subset {\mathbb Z}^{n_2}$.

\begin{proposition} \cite{H1,H3,H5}
	A~function $g:(X, k_1) \to (Y, k_2)$ is
	$(k_1, k_2)$-continuous if and only if for each point  $x\in X$,
	$g(N_{k_1}(x,1))\subset  N_{k_2}(g(x),1)$, where $X \subset {\mathbb Z}^{n_1}$ and $Y \subset {\mathbb Z}^{n_2}$.
\end{proposition}

Using  $k$-continuous maps,
we establish a {\it digital topological category}, denoted by {\it DTC},
consisting of the following two pieces of data \cite{H5}.\\
$\bullet$ The set of $(X, k)$ as  objects,\\
$\bullet$ For every ordered pair of objects  $(X, k_1)$ and $(Y, k_2)$,
the set of all $(k_1, k_2)$-continuous maps $f:(X, k_1) \to (Y, k_2)$ as morphisms.\\

Based on the digital continuity, since 2005 \cite{H4},
we have often used a {\it $(k_1, k_2)$-isomorphism} for a map between two sets $X \subset {\mathbb Z}^{n_1}$ and $Y \subset {\mathbb Z}^{n_2}$ in \cite{H4,HP1} instead of a {\it $(k_1, k_2)$-homeomorphism} for a map between two sets $X, Y \subset {\mathbb Z}^3$ in \cite{B2}, as follows:

\begin{definition} \cite{B2} (see also \cite{H4,HP1})
	Assume $(X, k_1)$ and $(Y, k_2)$ on ${\mathbb Z}^{n_1}$ and ${\mathbb Z}^{n_2}$, respectively. We say that
	a bijection $h: (X, k_1) \to (Y, k_2)$ is a
	$(k_1, k_2)$-isomorphism if $h$
	is a $(k_1, k_2)$-continuous and the inverse map of $h$,
	$h^{-1}: Y \to X$, is $(k_2, k_1)$-continuous.
	Then we denote this isomorphism by 
	$(X, k_1) \approx_{(k_1, k_2)}(Y, k_2)$.
	If $n_1 = n_2$ and $k_1 = k_2$, then we call it a {\it $k_1$-isomorphism}.
\end{definition}

Compared with a $(k_1, k_2)$-isomorphism,
a local $(k_1, k_2)$-isomorphism which is weaker than a  $(k_1, k_2)$-isomorphism was defined. 
The notion plays an important role in establishing a $(k_1, k_2)$-covering space \cite{H2,H3} and a pseudo-$(k_1, k_2)$-covering space \cite{H9,H17}, $DT$-$(k_1, k_2)$-embedding theorem \cite{H13}, and so forth.

\begin{definition} \cite{H1,H14}
	For two digital images $(X, k_1)$ on ${\mathbb Z}^{n_1}$ and $(Y, k_2)$ on ${\mathbb Z}^{n_2}$, consider a  map  $h: (X,k_1) \to  (Y,k_2)$. Then~the map $h$ is said to be a local (or $L$-, for brevity) $(k_1, k_2)$-isomorphism 
	if  for every $x \in X$, $h$ maps $N_{k_1}(x, 1)$ $(k_1, k_2)$-isomorphically onto $N_{k_2}(h(x), 1)$, i.e., the restriction map $h\vert_{N_{k_1}(x, 1)}:
	N_{k_1}(x, 1) \to N_{k_2}(h(x), 1)$ is a $(k_1, k_2)$-isomorphism.
	If $n_1 = n_2$ and $k_1 = k_2 $, then the map $h$ is called a local (or $L$-) $k_1$-isomorphism.
\end{definition}

In view of Definition 3, we observe that a local $(k_1, k_2)$-isomorphism is obviously a $(k_1, k_2)$-continuous map. However, the converse does not hold.
Compared with a local $(k_1, k_2)$-isomorphism, the following is defined.

\begin{definition} \cite{H2} (Radius $2$-local $(k_1,k_2)$-isomorphism)
	Assume a map  $h:(X, k_1) \to (Y, k_2)$, where $X \subset {\mathbb Z}^{n_1}$ and $Y \subset {\mathbb Z}^{n_2}$.
	Then we say that it is a radius $2$-local $(k_1,k_2)$-isomorphism if 
	for each point $x\in X$ the restriction $h$ to $N_{k_1}(x, 2)$, i.e., 
	$h\vert_{N_{k_1}(x, 2)}:N_{k_1}(x, 2) \to N_{k_2}(h(x), 2)$ is a 
	$(k_1,k_2)$-isomorphism.
	In particular, in the case of $k:=k_1=k_2$, we call it a radius $2$-local
	$k$-isomorphism.
\end{definition}

Comparing Definitions 3 and 4, we obtain the following:
\begin{remark} A radius $2$-local $(k_1,k_2)$-isomorphism is clearly a local $(k_1, k_2)$-isomorphism. However, the converse does not hold.
\end{remark}

\begin{example}  Consider the map $p:({\mathbb Z},2) \to SC_k^{n,l}:=(s_i)_{i \in [0,l-1]_{\mathbb Z}}$ given by $p(t)=s_{t(mod\,l)}$.\\
	(1) In the case
	of $l\geq 5$, the map $p$ is a
	radius $2$-local $(2,k)$-isomorphism.\\	
	(2) In the case of $l=4$, the map $p$ is not a radius $2$-local $(2,k)$-isomorphism but an $L$-$(2,k)$-isomorphism.	
\end{example}

A radius $2$-local $(k_1,k_2)$-isomorphism  strongly plays an important role in studying digital images from the viewpoint of digital homotopy theory and digital covering theory.
In detail, this notion is essential for establishing the homotopy lifting theorem \cite{H2} which can be used to calculate digital fundamental groups \cite{BK1,H2,H3}.\\
\begin{remark} (Utilities of  radius $2$-local $(k_1,k_2)$-isomorphism)\\
	$\bullet$ Establishment of the homotopy lifting theorem \cite{H2}\\
	$\bullet$ Calculation of digital fundamental groups \cite{BK1,H2,H3}
\end{remark}

\section{\bf Digital $k$-homotopy and digital $k$-homotopy equivalence}\label{s5}

Before 2005, the notion of a digital $k$-homotopy was considered in ${\mathbb Z}^3$. In detail, the first version of a digital $k$-homotopy was initially established by Khalimsky \cite{K1}.
After that, the second version of it was introduced by Boxer in ${\mathbb Z}^3$ \cite{B2}. 
After that, in ${\mathbb Z}^n, n \in {\mathbb N}$, the third version of it was 
established by Han in \cite{H2,H3}. Han's version is called a digital $k$-homotopy relative to a subset $A\subset X$ (see \cite{H2,H3}) and it has been used to study a $k$-homotopic thinning of a digital image $(X, k)$ in ${\mathbb Z}^n$.
For a digital image $(X, k)$ and $A \subset X$, we often call  $((X,A), k)$ a digital image pair. 
Let us now recall some history of them as follows:

\begin{definition} \cite{K1} (Khalimsky's version of a digital $k$-homotopy in ${\mathbb Z}^3$).	
	Let $(X, k_1)$ and $(Y, k_2)$ be digital images.
	Let $f, g: (X, k_1) \to (Y, k_2)$ be
	$(k_1,k_2)$-continuous functions.\\
	 Suppose that there exist $m \in {\mathbb N}$
	and a function  $H:X \times[0, m]_{\mathbb Z} \to Y$ such that\\
	(1) for all  $x\in X, H(x, 0)=f(x)$ and $H(x, m)=g(x)$ and\\
	(2) for all $t\in [0, m]_{\mathbb Z}$, the induced function $H_t: X \to Y$  given
	by
	$H_t(x) = H(x, t)$ for all $x \in X$ is $(k_1,k_2)$-continuous.\\
		Then we say that $H$ is a $k$-homotopy between $f$ and $g$.	
\end{definition}

As a more rigid version than the digital  $(k_1,k_2)$-homotopy of Definition 5, the following was constructed. In addition,
a pointed digital image is said to be a pair $(X, x_0)$
and $x_0 \in X$. A pointed $(k_1,k_2)$-continuous map
$f : (X, x_0) \to (Y, y_0)$ is a $(k_1,k_2)$-continuous function
from $X$ to $Y$ such that $f(x_0) = y_0$
from $(X, x_0)$ to $(Y, y_0)$ \cite{B2}. 

\begin{definition} \cite{B2} (Boxer's version of a (pointed) digital $(k_1,k_2)$-homotopy in ${\mathbb Z}^3$).	
		Let $(X, k_1)$ and $(Y, k_2)$ be a digital image, respectively.
	 Let $f, g: (X, k_1) \to (Y, k_2)$ be
	$(k_1,k_2)$-continuous functions. Suppose that there exist $m \in {\mathbb N}$
	and a function  $H:X \times[0, m]_{\mathbb Z} \to Y$ such that\\
	(1)-(2): the conditions (1) and (2) of Definition 5 and further, \\
	(3) for all $x \in X$, the induced function $H_x :[0, m]_{\mathbb Z} \to Y$
	given by $H_x(t)= H(x, t)$ for all $t \in [0,m]_{\mathbb Z}$ is $(2, k_2)$-continuous.\\ 
	Then we say that $H$ is a $(k_1,k_2)$-homotopy between $f$ and $g$.\\
A digital $(k_1,k_2)$-homotopy
$H : X \times [0,m]_Z \to Y$
between $f$ and $g$ is called a pointed digital $(k_1,k_2)$-homotopy
between $f$ and $g$ if for all $t \in [0,m]_Z, H(x_0, t) = y_0$.
		\end{definition}
	
As a generalization of the Boxer's pointed digital $(k_1,k_2)$-homotopy, Han \cite{H2,H3} established the so-called $\lq\lq$digital $(k_1, k_2)$-homotopy relative to $A$" in \cite{H2,H3} with the expansion of the dimension from ${\mathbb Z}^3$ to ${\mathbb Z}^n, n \in {\mathbb N}$ and $\{x_0\}$ to $A\subset X$, as follows:

\begin{definition} \cite{H1,H5} (Han's version of a digital $(k_1, k_2)$-homotopy relative to $A$ in ${\mathbb Z}^n, n\in {\mathbb N}$).	
	Let $(X, k_1)$ and $(Y, k_2)$ be digital images on ${\mathbb Z}^{n_1}$ and ${\mathbb Z}^{n_2}$, respectively. 
	 Let $f, g: (X,k_1) \to (Y, k_2)$ be 
	$(k_1, k_2)$-continuous functions. Suppose that there exist $m \in {\mathbb N}$
	and a function  $H:X \times[0, m]_{\mathbb Z} \to Y$ such that\\
	(1)-(3): the conditions (1)-(3) of Definition 6, and further,\\
	(4) for all $t\in [0, m]_{\mathbb Z}$, assume that the induced map $H_t$ on $A$ is a constant
	which follows the prescribed function from $A$ to $Y$. To be precise, $H_t(x)=f(x)=g(x)$
	for all $x\in A$ and for all $t\in [0, m]_{\mathbb Z}$.\\
	Then we call $H$ a $(k_1, k_2)$-homotopy relative to $A$ between $f$ and $g$, and  we say that $f$ and $g$ are
	$(k_1, k_2)$-homotopic relative to $A$ in $Y$, $f\simeq_{(k_1, k_2)rel A}g$ in symbols.\\
	In particular, in the case of $n_1=n_2$ and $k:=k_1=k_2$, we abbreviate  $\lq\lq$$(k_1, k_2)$-homotopy relative to $A$" as $\lq\lq$$k$-homotopy relative to $A$" and use the notation  $f\simeq_{k rel A} g$.
\end{definition}

In Definition 7, if $A= \{x_0\}\subset X$, then we observe that
$H$ is a pointed $(k_1, k_2)$-homotopy at $\{x_0\}$ \cite{B2}. If, for some $x_0\in X$, $1_X$ is $k$-homotopic to the constant map in the space $x_0$ relative to $\{x_0\}$, then
we say that $(X, x_0)$ is pointed $k$-contractible \cite{B2}.
For instance, $SC_{3n-1}^{n, 4}$ is $(3n-1)$-contractible \cite{B2,H3,H6}.

\begin{remark} Using Han's homotopy of Definition 7, we can proceed a $k$-homotopic thinning \cite{H6} (see Definition 9 below). 
\end{remark}

In \cite{H6}, we have studied a $k$-homotopic
thinning in $DTC$ by
using a strong $k$-deformation
retract derived from the $k$-homotopy relative to some subset $(A, k) \subset (X, k)$ in Definition 7.

\begin{definition} \cite{H6} In DTC, for a digital image pair $((X, A), k)$,
	$A$ is said to be a strong $k$-deformation retract of $X$ if there
	is a $k$-retraction $r$ of $X$ onto $A$ such that $F:i \circ r
	\simeq_{k \cdot {\it rel.}A} 1_X$, where $i:(A,k) \to (X,k)$ is an inclusion map.
	Then, a point $x \in X \setminus A$ is called strongly $k$-deformation retractable.
\end{definition}

\begin{definition} \cite{H6} In $DTC$, for a digital image $(X, k)$, a process of deleting strong $k$-deformation retractable points
	from $X$ is called a $k$-homotopic thinning.
\end{definition}

Based on  Definition 7, the notion of a digital $(k_1, k_2)$-homotopy equivalence is introduced as follow:

\begin{definition} \cite{HP2} 
Assume two digital images $(X, k_1)$ on ${\mathbb Z}^{n_1}$ and $(Y, k_2)$ on  ${\mathbb Z}^{n_2}$. 
	Assume a $(k_1, k_2)$-homotopy from $(X, k_1)$ and $(Y, k_2)$ such that 
	$g\circ f \simeq_{(k_1,k_2))} 1_{X}$ and a $(k_2, k_1)$-homotopy from $(Y, k_1)$ and $(X, k_1)$ such that 
	$f\circ g \simeq_{(k_2,k_1))} 1_{Y}$.
	Then we say that $(X, k_1)$ is $(k_1, k_2)$-homotopy equivalent to  $(Y, k_2)$.
	In particular, in the case of $n_1=n_2$ and $k:=k_1=k_2$, the $(k_1, k_2)$-homotopy equivalence is said to be a $k$-homotopy equivalence.
	\end{definition}

Even though Boxer often claimed that the term $\lq\lq$$k$-homotopy equivalence" was introduced in \cite{B2}, one can clearly observe that there is no notion (or term) of a $k$-homotopy equivalence in \cite{B2}. The paper \cite{B2} just only stated an equivalence relation of a homotopy of Definition 6 (see Proposition 2.8 of \cite{B2}).
As for a digital $(k_1,k_2)$-homotopy equivalence, Han \cite{HP2} delivered his talk on 
a digital $(k_1,k_2)$-homotopy equivalence and a digital $k$-homotopy equivalence at the international conference which was held in Washington in the U.S. in 2003 (see \cite{HP2}).
Besides, many scholars including the professors R. Kopperman, J. D. Lawson, M. Mislove enjoyed the talk and discussed the notion of a digital $(k_1,k_2)$-homotopy equivalence and a digital $k$-homotopy equivalence there.
After two years later, in only 2005, Boxer studied the notion of a digital $(k_1,k_2)$-homotopy equivalence and a digital $k$-homotopy equivalence  in his paper.\\

\begin{remark}  Based on Khalimsky topology or an adjacency derived from the Khalilmsky topology on ${\mathbb Z}^n$, the paper \cite{H8} studied many types of homotopy equivalence from the viewpoint of digital topology.
	Besides, the paper \cite{HY1} also studied a special kind of homotopy equivalence using some adjacency derived from the Marcus-Wyse topology.
\end{remark}

 \section{\bf Further remarks}\label{s8}
 
 We have discussed some utilities of the notion of a radius $2$-local $(k_1,k_2)$-isomorphism and referred to a history of the concept of a digital $k$-homotopy.
 Based on the works \cite{B2,B3,BK1,H1,H2,H3,H5,H6,H14,H17,HP2,K1}, we will remark on Boxer's misconception in \cite{B3} on Han's papers \cite{H9,H14} and a revision of a pseudo-$(k_1, k_2)$-covering space. 
 In particular, since Boxer \cite{B3} referred to some results in \cite{H17} in an inaccurate way, we will correct some assertions in \cite{B3} which are due to Boxer's misunderstanding with a next paper as a series of this topic.\\

 {\bf Conflicts of Interest}: The author declares no conflict of interest.

\end{document}